\documentclass{siamltex}

\usepackage{amsmath}
\usepackage{graphicx}
\usepackage{amssymb}
\usepackage{dsfont}
\usepackage{wrapfig}
\usepackage{color}
\usepackage{multirow}
\usepackage{subfigure}
\usepackage{url}
\usepackage{float}
\usepackage{bm}
\usepackage{lineno}

\usepackage[
colorlinks=true,
citecolor=blue,
urlcolor=blue,
linktocpage,
pagebackref
]{hyperref}

\allowdisplaybreaks
\numberwithin{equation}{section} \numberwithin{figure}{section}
\numberwithin{table}{section}

\newtheorem{example}{Example}[section]

\renewcommand{\div}{\operatorname{div}}

\newcommand{\Ker}{\operatorname{Null}}

\newcommand{\Ran}{\operatorname{Range}}

\newcommand{\wh}[1]{\widehat{#1}}

\newlength\figureheight
\newlength\figurewidth
\newcommand{\curl}{\operatorname{curl}}
\renewcommand{\div}{\operatorname{div}}

\def\b1{{\mathbf 1}}
\def\bv{{\mathbf v}}
\def\bu{{\mathbf u}}

\def\bg{{\mathbf g}}

\def\bx{{\mathbf x}}
\def\by{{\mathbf y}}

\def\bn{{\mathbf n}}
\def\bbf{{\mathbf f}}

\def\b0{{\mathbf 0}}

\renewcommand{\O}{{\mathcal O}}

\def\T{{\mathcal T}}

\def\l1{{$\ell_1$}}

\def\hA{\widehat{A}}
\def\hB{\widehat{B}}

\def\blambda{{\boldsymbol \lambda}}

\def\bone{{\boldsymbol 1}}
\newcommand{\vertiii}[1]{{\left\vert\kern-0.25ex\left\vert\kern-0.25ex\left\vert #1
    \right\vert\kern-0.25ex\right\vert\kern-0.25ex\right\vert}}
\newcommand{\jump}[1]{\left [ #1 \right ]}
\def\hM{\widehat M}
\def\hW{\widehat W}
\def\hx{\widehat \bx}
\def\hy{\widehat \by}

\def\hY{\widehat Y}
\def\giter{GMRES iterations}
\def\fires{Final residual $\times$ 1e-5}
\def\hf{\widehat \bbf}

\definecolor{tk}{rgb}{0.77,0,0.04}    

\definecolor{pv}{rgb}{0,0.52,0.82}   

\begin{document}

\title{Algebraic Hybridization and Static Condensation with Application to
Scalable $H(\div)$ Preconditioning
\thanks{
  This work was performed under the auspices of the U.S. Department of Energy by
  Lawrence Livermore National Laboratory under Contract DE-AC52-07NA27344, and
  sponsored by the U.S. Department of Energy, Office of Science, Office of
  Advanced Scientific Computing Research, Applied Mathematics
  program. LLNL-JRNL-732140.}
}


\author{
V.~Dobrev%
\thanks{Center for Applied Scientific Computing,
  Lawrence Livermore National Laboratory,
  P.O. Box 808, L-561,
  Livermore, CA 94551, U.S.A (dobrev1@llnl.gov,
  tzanio@llnl.gov, cslee@llnl.gov, tomov2@llnl.gov, panayot@llnl.gov).}
\and T.~Kolev%
\footnotemark[2]
\and C.~S.~Lee%
\footnotemark[2]
\and V.~Tomov%
\footnotemark[2]
\and P.~S.~Vassilevski%
\footnotemark[2]
}

\maketitle

\begin{abstract}
We propose an unified algebraic approach for static condensation and
hybridization, two popular techniques in finite element discretizations.  The
algebraic approach is supported by the construction of scalable solvers for
problems involving $H(\div)$-spaces discretized by conforming (Raviart-Thomas)
elements of arbitrary order. We illustrate through numerical experiments the
relative performance of the two (in some sense dual) techniques in comparison
with a state-of-the-art parallel solver, ADS \cite{kv12}, available in
\cite{hypre} and \cite{mfem}.  The superior performance of the hybridization
technique is clearly demonstrated with increased benefit for higher order
elements.
\end{abstract}

\begin{keywords}
static condensation, hybridization, algebraic multigrid, mixed finite elements,
ADS, $H(\div)$ solvers, radiation-diffusion transport.
\end{keywords}

\begin{AMS}
65F10, 65N20, 65N30
\end{AMS}

\pagestyle{myheadings} \thispagestyle{plain}
\markboth{DOBREV, KOLEV, LEE, TOMOV, and VASSILEVSKI}
{{Algebraic hybridization and static condensation}}

\section{Introduction}\label{section: introduction}

This paper proposes a unified algebraic approach of constructing reduced
problems by two popular techniques in the finite element method: static
condensation and hybridization. It also discusses approaches for constructing
scalable preconditioners for these reduced problems. The present work is
motivated by the solution of problems arising in variety of finite element
simulations and of particular interest is the case when the problems involve the
function space $H(\div)$ ($L_2$-vector functions that have their weak divergence
in $L_2$).  Such problems naturally arise in the mixed finite element method for
Darcy equation \cite{BoffiBrezziFortin2013}, Brinkman equation \cite{vv13}, as
well as in the FOSLS (first order system least-squares) discretization approach,
\cite{clmm94}.  Also, certain formulations of radiation-diffusion transport
\cite{Brunner02}, lead to problems involving the space $H(\div)$. This is in
fact one of the practical examples from a more realistic simulation that we
consider in our numerical tests.

Highly scalable solvers for $H(\div)$ problems have been designed in the past.
The most successful ones, of ADS-type \cite{kv12}, were based on the
{\em regular decomposition} theory developed by Hiptmair and Xu \cite{hx07}.
Although of optimal complexity, the current state-of-the-art ADS algorithms are
quite involved and require also going through solvers for $H(\curl)$, such as
AMS, \cite{kolev-vassilevski-ams}. Both solvers, ADS and AMS, require some
additional user input (e.g., discrete gradient), which makes them not completely
{\em algebraic}.  The results of the present paper can be viewed as a further
step towards the design of fully algebraic solvers for problems involving
$H(\div)$.

To design faster, yet still scalable, alternative to ADS, we employ a
traditional finite element technique commonly used in mixed finite elements
referred to as {\em hybridization}. This is perhaps the first technique proposed
to solve the saddle-point problems arising in the mixed methods, since it leads
to symmetric and positive definite (SPD) reduced problems. Its early form
appeared in a paper by Fraeijis de Veubekein in 1965 \cite{veubeke65} as an
efficient solution strategy for elasticity problems. In \cite{ArnoldBrezzi85},
Arnold and Brezzi studied the method from a theoretical perspective, and
obtained error estimates for all the unknowns involved; see also
\cite{bm94}. The technique was further developed by Cockburn and Gopalakrishnan
to allow approximations by polynomials of different degree in different
subregions; see for example \cite{cg05} and \cite{cg05_GAMM} for applications to
mixed finite element approximations of second order scalar elliptic problems and
Stokes flow.

Another very popular technique in finite elements is {\em static condensation}.
It is a reduction technique intended to reduce the size of a linear system of
finite element equations by eliminating at element or subdomain level internal
degrees of freedom (dofs). It was also first introduced in the structural
analysis literature in 1965, cf., \cite{guyan65, irons65, wilson65}, and has
been widely used since then. Static condensation neglects the dynamic effect in
the reduction, hence the name. On element level, one can eliminate locally
degrees of freedom that are not shared by neighboring elements (a typical case
for high order elements). Thus, the reduced, Schur complement, problem contains
only shared (interface) dofs. If the elimination is performed on a subdomain
level, the resulting reduced system coincides with the one in the balancing
domain decomposition by constraints (BDDC) algorithm \cite{d03}. Recently, BDDC
preconditioners with deluxe scaling and adaptive selection of primal constraints
for $H(\div)$ problems have been developed, see for instance \cite{Oh17,
Zampini16, ZampiniKeyes16}. The approach in the present paper, on the other
hand, solves the reduced system by algebraic multigrid (AMG).

For $H(\div)$-conforming elements (i.e., elements that have only shared dofs
through common faces), both static condensation and hybridization lead to
reduced problems of the same size (equal to the number of shared interface
dofs). If one applies these two approaches not on element level, but for
subdomains (union of elements), substantial reduction in size can be
achieved. Thus, these two approaches are of great practical interest and the
goal of the present paper is to study them in a common framework, emphasize
their similarities and differences, and most importantly, design new, or modify
existing solution techniques that are efficient and scalable, achieving
substantial savings compared to the state-of-the-art solvers for the original
problems. For problems obtained by static condensation, there has been success
in modifying the state-of-the-art solvers (AMS and ADS) to be directly
applicable to the reduced problem, cf., \cite{brunner-kolev-amg-explicit,
dpg-2016}. In this paper, we focus on the design of scalable solvers for
problems involving $H(\div)$ in a general algebraic setting, extending the
preliminary results in \cite{LeeVassilevski17}.

The remainder of the paper is organized as follows. In Section
\ref{section:hybridiztion}, we describe the algebraic hybridization approach and
show how it can serve as an alternative to traditional finite element assembly.
Next, we deduce its relation to static condensation and the mixed hybridized
finite element method in Section \ref{section: static condensation}.
Preconditioning strategies for the reduced systems obtained from hybridization
and static condensation are discussed in Section \ref{section:preconditioning}.
Then we describe the implementation details of the methods and present some
numerical examples in Section \ref{section: implementation} and Section
\ref{section: numerical tests} respectively. Finally, we conclude the paper in
Section \ref{section: conclusion}.


\section{Two perspectives of finite element assembly}
\label{section:hybridiztion}

In this section we describe the traditional process of finite element assembly,
and introduce algebraic hybridization as an alternative process to obtain a
global linear system describing the same finite element discretization.

\subsection{Finite element assembly components}
\label{subsection: components of finite element assembly}

To fix notation, consider a bilinear form $a(\cdot,\cdot)$ and a linear form
$f(.)$ over a Hilbert space $\mathcal{H}$ defining the Galerkin weak variational
formulation: Find $u \in \mathcal{H}$, such that for all $v \in \mathcal{H}$, we have
\begin{equation*}
a(u,v) = f(v).
\end{equation*}

In the finite element method, we introduce a triangulation $\T=\T_h$ of the
given computational domain and decompose the forms into sums of their local
element versions, i.e $a(\cdot,\cdot) = \sum\limits_{\tau \in \T}
a_\tau(\cdot,\cdot)$ and $f(.) = \sum\limits_{\tau \in \T} f_\tau(.)$.  We
further replace $\mathcal{H}$ with a finite dimension space induced by $\T_h$ and on each
element $\tau \in \T$ we compute a local {\em element stiffness matrix} $A_\tau$
and a {\em load vector} $\bbf_\tau$ by using the element basis functions in the
local bilinear and linear forms.  We are assuming that $A_\tau$ (which
correspond to $a_\tau(\cdot,\cdot)$) are SPD. The latter holds for example if
the bilinear form $a(u,v)$ contains a mass term, $\int \beta uv\;dx$, with a
positive coefficient $\beta=\beta(x) > 0$.

\subsection{The traditional finite element assembly procedure}
\label{subsection: traditiona lassembly}

A central concept in the finite element method is that of {\em assembly} -- the
accumulation of local stiffness matrices $A_\tau$ and load vectors $\bbf_\tau$
into a global matrix $A$ and a global vector $\bbf$ which form the (global)
linear system
\begin{equation}\label{Axf}
A \bx = \bbf,
\end{equation}
for the unknown vector of finite element degrees of freedom.  The $n \times n$
matrix $A$ is SPD and generally sparse, i.e., the number of its nonzero entries,
$nnz(A)$, is $\O(n)$.

Let ${\widehat A} = blockdiag(A_\tau)$ be the $\hat{n} \times \hat{n}$
block-diagonal matrix with all element stiffness matrices on its main diagonal,
and $\widehat \bbf$ be similarly the $\hat{n} \times 1$ vector containing all
element load vectors.  The finite element method produces a global-to-local
mapping, $P$, which maps global dofs into local ones for each element.  That is,
for a $n \times \hat{n}$ matrix $P$, the traditional assembly process can be
written as the triple matrix product that relates $A$ and ${\widehat A}$, and a
matrix-vector product that builds the global right-hand side,
\begin{equation}\label{assembly}
A = P^T \widehat A P\,,
\quad\text{and}\quad
\bbf = P^T \widehat \bbf \,.
\end{equation}
In the simplest case the action $\widehat \bx = P \bx$ simply means copying the
values of the global dof vector $\bx$ to the corresponding local values,
independently in each element. More complicated $P$ can account for
non-conforming mesh refinement as we describe in some detail in a following
section.  Either way, \eqref{assembly} is a formalized expression for the
traditional process of finite element assembly.

\subsection{Algebraic hybridization}
\label{subsectioN: algebraic hybridization}

It is classical in the finite element literature to introduce \eqref{Axf} as the
global linear system that corresponds to the finite element discretization
specified by ${\widehat A}$, $\widehat \bbf$ and $P$.  As we describe next,
there is an alternative global linear system corresponding to the same solution,
which in certain cases is advantageous from the perspective of linear solvers.
To introduce this alternative, let $C$ be another matrix such that
\begin{equation*}
\Ker(C) = \Ran(P),
\end{equation*}
i.e. a vector $z$ satisfies $C z = 0$, if and only if there is a vector $y$,
such that $z = P y$.  In particular, we have the matrix equality $C P = 0$.

We assume that $P$ is full column rank and that $C$ is full row rank. Similarly
to $P$, the matrix $C$ could also be provided algebraically by the finite
element space, e.g. the columns of $C^T$ can be constructed by orthogonally
completing the basis given by the columns of $P$. (The matrix $C$ is not unique
and there are other ways to construct it as illustrated later in the paper.)

The idea of algebraic hybridization is that instead of adding the different
element contributions of ${\widehat A}$, $\widehat \bbf$, we use $C$ to enforce
equality constraints of the decoupled vector $\widehat \bx$ through a new
Lagrange multiplier variable $\blambda$.  This process is known as {\em
  hybridization} in the context of mixed finite element method \cite{bf91}, but
in this section we consider it on a purely algebraic level, independent of the
particular continuous problem or finite element discretization approach.  In
fact, the $C$-based reformulation can be described in more general terms than
classical assembly, which we do in the following main linear algebra result.

\begin{theorem}\label{theorem: equivalence of solutions}
Let ${\widehat A}$, $P$ and $C$ be matrices of size $\hat{n} \times
\hat{n}$, $\hat{n} \times n$ and $m \times
\hat{n}$ ($m=\hat{n}-n$) respectively, such that ${\widehat A}$
is SPD, $P$ is full column rank, $C$ is full row rank, and $\Ker(C) = \Ran(P)$.
Also, let $\widehat \bbf$ be a column vector of size $\hat{n}$. Then the problem
\begin{equation}\label{factored original system}
P^T{\widehat A} P\, \bx = P^T\, {\widehat \bbf},
\end{equation}
is equivalent to the following, {\em hybridized}, saddle-point system
\begin{equation}\label{hybridized system}
\begin{bmatrix}
{\widehat A} & C^T \\
C & 0
\end{bmatrix}
\begin{bmatrix}
{\widehat \bx}\\
\blambda
\end{bmatrix} =
\begin{bmatrix}
{\widehat \bbf}\\
0
\end{bmatrix}.
\end{equation}
Specifically, the saddle-point system \eqref{hybridized system} is
uniquely solvable and the solution ${\widehat \bx}$ of
\eqref{hybridized system} and $\bx$ (the solution of the original
system), are related as
\begin{equation*}
{\widehat \bx} = P \bx.
\end{equation*}
The solution $\widehat \bx$ is computable via standard (block-)Gaussian elimination:
we first compute the Schur complement system for the Lagrange multiplier
\begin{equation}\label{Lagrange multiplier system}
H \blambda \equiv C \wh{A}^{-1} C^T \blambda = C \wh{A}^{-1} \wh{\bbf},
\end{equation}
then, by back substitution, we have
\begin{equation*}
{\widehat \bx} = {\widehat A}^{-1}({\widehat \bbf}-C^T\blambda).
\end{equation*}
\end{theorem}
\begin{proof}
The saddle-point system is invertible due to our assumptions; namely,
${\widehat A}$ is SPD (hence invertible), and $C^T$ is full column
rank, hence the (negative)
Schur complement $H = C {\widehat A}^{-1} C^T$ is also
SPD and invertible.
From ${\widehat A} {\widehat \bx} + C^T \blambda = {\widehat \bbf}$ and
$CP =0$, we obtain
\begin{equation*}
P^T {\widehat A} {\widehat \bx} = P^T {\widehat \bbf} =: \bbf.
\end{equation*}
Now, since $C {\widehat \bx} = 0$ (the second equation of
\eqref{hybridized system}), and from the assumption $\Ker(C) = \Ran(P)$,
it follows that there is an $\bx$ such that ${\widehat \bx} = P \bx$.
Then from $\bbf = P^T {\widehat A} {\widehat \bx} = P^T {\widehat A}
P\bx$, we see that this $\bx$ is unique; namely, $\bx$ is the unique
solution of the original problem \eqref{factored original system}.
\end{proof}

We can summarize the results of the theorem and the previous section by stating
that $P^T{\widehat A} P$ and $C {\widehat A}^{-1} C^T$ are equally valid global
matrices for the same problem.  In the case of finite elements, the first one
corresponds to assembly, while the second one is the hybridized system for
Lagrange multipliers.

Note that the condition that $P$ has full column rank is equivalent to the
existence of a left inverse of $P$: there is an $n \times \hat{n}$ matrix $R$
such that $R P = I$, e.g.\ $R=(P^T P)^{-1} P^T$.  We can use this restriction
matrix $R$ to compute the solution $\bx$ of \eqref{factored original system}
from the solution $\widehat\bx$ of \eqref{hybridized system} as $\bx = R
\widehat\bx$. Furthermore, if instead of $\widehat\bbf$ we have access only to
the assembled right-hand side vector $\bbf = P^T \widehat\bbf$, we can use
$\widetilde\bbf = R^T \bbf$ in \eqref{hybridized system} instead of
$\widehat\bbf$ since $P^T \widetilde\bbf = P^T \widehat\bbf = \bbf$.

Note also that the first set of equations in \eqref{hybridized system},
$$
{\widehat A} {\widehat \bx} = {\widehat \bbf} - C^T\blambda.
$$
can be interpreted as a local version of the problem (in a weak form), where the
term $C^T\blambda$ plays the role of a dual vector for Neumann-type boundary
conditions.
Indeed, if we integrate by parts locally, we get
$$
{\widehat A} {\widehat \bx} = {\widehat \bbf} - {\widehat F}
$$
where ${\widehat F}$ is the Neumann data for ${\widehat \bx}$ tested against
the test functions in each element. Since the exact solution satisfies
\eqref{factored original system}, we have $P^T {\widehat F} = 0$, i.e.
$\widehat F = C^T \blambda$ for some $\blambda$ by $\Ker(P^T)=\Ran(C^T)$.
One can then think of hybridization as prescribing Neumann data for local
problems from one set of $\blambda$ values and then imposing the fact that the
local solutions (fluxes) should match on their shared dofs to derive equations
for $\blambda$; for more details, see \cite{misha}.

The result of Theorem~\ref{theorem: equivalence of solutions} can be extended to
more general settings that relax the conditions that $\widehat{A}$ and $P$ are
full column rank. This is illustrated in the next theorem.

\begin{theorem}\label{theorem: general equivalence of solutions}
Assume that $\widehat{A}$ is a square matrix (without assuming that it is SPD or
even invertible) and that $\Ker(C)=\Ran(P)$ (without assuming that $P$ has full
column rank and $C$ has full row rank). Then equations
\eqref{factored original system} and \eqref{hybridized system} are equivalent in
the following sense
\begin{itemize}
\item
If $\bx$ solves \eqref{factored original system} then there exists $\blambda$
such that $(\widehat\bx=P\bx,\blambda)$ solves \eqref{hybridized system}.
\item
If $(\widehat\bx,\blambda)$ solves \eqref{hybridized system} then there exists
$\bx$ that solves $\eqref{factored original system}$ such that
$P\bx = \widehat\bx$.
\end{itemize}
If, in addition, we assume that $\widehat{A}$ is invertible then equations
\eqref{hybridized system} and \eqref{Lagrange multiplier system} are also
equivalent in the following sense
\begin{itemize}
\item
If $(\widehat\bx,\blambda)$ solves \eqref{hybridized system} then $\blambda$
solves \eqref{Lagrange multiplier system}.
\item
If $\blambda$ solves \eqref{Lagrange multiplier system} then
$(\widehat\bx = \widehat{A}^{-1}(\widehat\bbf - C^T\blambda),\blambda)$ solves
\eqref{hybridized system}.
\end{itemize}
\end{theorem}
\begin{proof}
Assume that $\bx$ solves \eqref{factored original system} and set
$\widehat\bx=P\bx$. We need to show that there exists $\blambda$ such that
$(\widehat\bx,\blambda)$ solves \eqref{hybridized system}. The second equation
in \eqref{hybridized system} follows from $C \widehat\bx = C P \bx = 0$. The
defining condition for $\blambda$ is
\[
\widehat{A} \widehat\bx + C^T \blambda = \widehat\bbf \,,
 \qquad\text{or}\qquad
C^T \blambda = \widehat\bbf - \widehat{A} \widehat\bx =: \widehat\bg \,.
\]
By \eqref{factored original system}, we have that
$P^T \widehat\bg = 0$, i.e.\ $\widehat\bg\in\Ker(P^T)$. Since we have
\[
\Ker(P^T) = \Ran(P)^\perp \subseteq \Ker(C)^\perp = \Ran(C^T) \,,
\]
the existence of the required Lagrange multiplier $\blambda$ is proven.

Assume now that $(\widehat\bx,\blambda)$ solves \eqref{hybridized system}.
On one hand, the second equation in \eqref{hybridized system} is
$C \widehat\bx=0$, i.e.\ $\widehat\bx\in\Ker(C)$ and since
$\Ker(C)\subseteq \Ran(P)$, there exists $\bx$ such that $P\bx=\widehat\bx$.
On the other hand, the first equation in \eqref{hybridized system}, implies
\[
P^T (\widehat{A} \widehat\bx + C^T \blambda ) = P^T \widehat\bbf \,,
  \qquad\text{or}\qquad
P^T \widehat{A} \widehat\bx = P^T \widehat\bbf \,,
\]
showing that $\bx$ solves \eqref{factored original system}.

The proof of the second part of the theorem is straightforward.
\end{proof}

Note that the uniqueness of $\bx$ is not guaranteed in the settings of Theorem
\ref{theorem: general equivalence of solutions}, even if we assume that $P$ has
full column rank and $\widehat{A}$ is invertible. Similarly, the uniqueness of
$\blambda$ does not follow from the additional assumption that $C$ has full row
rank and $\widehat{A}$ is invertible.

\subsection{A main example of matrices $P$ and $C$}

Below we discuss an application of the main theorem to a slightly more general
case of finite element discretization.  This is the main example we have in mind
in the application to $H(\div)$ bilinear form considered in the present paper.

\begin{example}\label{example: P and C}
In a finite element setting, we may think of the following $P$ originating from
a generally non-conforming finite element space.

Assume, that we have elements (or more generally subdomains that are unions of
finite elements).  We also have interfaces between these substructures (original
elements or subdomains). There are degrees of freedom interior, ``$i$'', to the
substructures. We also have dofs on the interfaces. In this setting we assume
that on the interfaces we have dofs that are {\em master}, ``$m$'', and {\em
slave}, ``$s$''.

Introduce the matrices $P$ and $C$ with the described blocks of $i, m,$ and
$s$-dofs,
\begin{equation*}
P =
\begin{bmatrix}
 I_i & 0 \\
  0  & I_m \\
 0   & W_{sm}
\end{bmatrix}
\qquad \text{and} \qquad
C = \begin{bmatrix} 0 & -W_{sm} & I_s \end{bmatrix}.
\end{equation*}
The 3-by-2 block-form of $P$ above, corresponds to the interior and master dofs
(for its columns) and interior, master and slave dofs (for its rows). The block
$W$ represents the mapping that determines the slave dofs in terms of the master
ones. Finally, the matrix ${\widehat A}$, is the block-diagonal matrix, with
blocks assembled for the individual substructures (or simply the element
matrices).

As before, we assume that ${\widehat A}$ is invertible, which is the case if we
have a mass term added to a semidefinite term, like the $H(\div)$ bilinear form
we consider in this paper.
\end{example}

\begin{example}\label{example: C from bilinear form}
Consider the case when the solution space is $H(\div)$-conforming, and introduce
a multiplier finite element space, $\mathcal{M} \subset L^2(\mathcal{I})$ on the
interface $\mathcal{I}$ between elements (or substructures), see Figure
\ref{fig:rt1-example} for illustration.

If we denote the {\em broken} solution space by $\wh{\mathcal{V}}$, then the
matrix $C$ can be derived from the bilinear form:
\[
c(\hat\bu,\lambda) =
  \int\limits_{\mathcal{I}} \jump{\hat\bu \cdot \bn} \lambda \,,\qquad
  \hat\bu \in \wh{\mathcal{V}}, \lambda\in\mathcal{M} \,.
\]
Above, $ \jump{\hat\bu \cdot \bn}$ stands for the jump of the normal component
of $\hat\bu$ across the interface ${\mathcal I}$.

Similarly to Example \ref{example: P and C}, for the same $P$, the matrix $C$
now has the $M_s$-weighted form
\[
C = \begin{bmatrix} 0  &  -M_s W_{sm} & M_s \end{bmatrix} \,,
\]
where $M_s$ is a mass matrix between the multiplier space $\mathcal{M}$ and the
normal trace space of $\wh{\mathcal{V}}$.
\end{example}

\section{Relation to static condensation and hybridized mixed finite elements}
\label{section: static condensation}

In this section we show that the hybridized matrix $H$, \eqref{Lagrange
multiplier system}, is closely related to the matrices obtained with two other
classical finite element approaches: that of algebraic {\em static condensation}
and that of traditional hybridization used in the mixed finite element method
for second order elliptic equations.

\subsection{Static condensation}\label{subsectioN: static condensation}

In the setting of Example \ref{example: P and C}, the matrix ${\widehat A}$ can
be partitioned into a two-by-two block from with its first block corresponding
to the interior dofs (with respect to the substructures). The second block
corresponds to the master and slave dofs combined, i.e., we have
\begin{equation} \label{A-2x2}
{\widehat A} = \left [
\begin{array}{cc}
{\widehat A}_{ii} & {\widehat A}_{ib}\\
{\widehat A}_{bi} & {\widehat A}_{bb}
\end{array} \right ].
\end{equation}
We notice that ${\widehat A}_{ii}$ is block-diagonal, since the interior dofs in
one substructure do not couple with interior dofs in any other substructure. The
Schur complement
\begin{equation*}
{\widehat S} = {\widehat A}_{bb} -
    {\widehat A}_{bi} {\widehat A}^{-1}_{ii} {\widehat A}_{ib},
\end{equation*}
is SPD, hence invertible. It acts on vectors corresponding to the interface dofs
``$b$'' (combined slave and master ones).
Now, consider the Schur complement $H = C {\widehat A}^{-1} C^T$
from \eqref{Lagrange multiplier system} for the hybridized system
\eqref{hybridized system}.
Since $C = [0, \;-W,\; I]$, and
\begin{equation*}
{\widehat A}^{-1} = \left [
\begin{array}{cc}
* & * \\
* & {\widehat S}^{-1}
\end{array} \right ],
\end{equation*}
we easily see that
\begin{equation*}
H = C {\widehat A}^{-1}C^T = [-W,\;I] {\widehat S}^{-1} [-W,\;I]^T,
\end{equation*}
while the {\em static condensation} matrix, i.e., the reduced matrix for the
master dofs is
\begin{equation}\label{static condensation matrix}
S = \left [
\begin{array}{c}
I\\
W
\end{array} \right ]^T{\widehat S} \left [
\begin{array}{c}
I\\
W
\end{array} \right ].
\end{equation}
In the case of finite element matrices, when the coupling across elements (or
more generally, across substructures) occurs only via interface dofs (no vertex
and edge dofs), one can decouple each interface dof into exactly two copies, and
let $W = -I$.  Examples of such elements are the Raviart--Thomas
($H(\div)$-conforming elements); and also, the $H^1$-non-conforming
Crouzeix-Raviart elements.

In this setting ${\widehat A}$ is block-diagonal, and also ${\widehat S}$ is
block-diagonal (substructure-by-substructure).  Thus, to form the hybridized
Schur complement $H$, we need to assemble the inverses of the local substructure
Schur complements, whereas to form the static condensation matrix $S$, we need
to assemble the local substructure Schur complements.

In other words, we may view the hybridization approach as a sort of a ``dual''
technique to the static condensation approach.  This has important consequences
for the construction of solvers for the respective problems. For example, we can
expect that solvers for $S$ will have the same nature as solvers for the
original matrix $A$, whereas solvers for $H$, will need preconditioners that are
effective in the dual to the original trace space.  This is the case for the
$H(\div)$ problem, as observed in \cite{LeeVassilevski17}.

\subsection{Hybridized mixed finite element method}\label{subsection: hybridized mixed f.e. method}

We next consider traditional hybridization for mixed finite elements,
cf. \cite{LeeVassilevski17}, applied in the settings of the $H(\div)$ problem:
find $\bu\in H(\div;\Omega)$ such that
\begin{equation}
a_{\div}(\bu,\bv) = (\alpha\, \div\bu, \div\bv) + (\beta \bu, \bv) = (\bg, \bv)
    \quad\quad\forall\;\bv\in H(\div;\Omega),
\label{eq:weak_hdiv}
\end{equation}
where the coefficients $\alpha$ and $\beta$ are positive and in general
heterogeneous.  In this case, $\hA$ has the following special form
\begin{equation}
\hA = \hB^T \hW^{-1} \hB + \hM.
\label{eq:discrete_hdiv}
\end{equation}
By introducing an additional unknown $q=\div\bu$, a mixed formulation of
\eqref{eq:weak_hdiv} can be written as: find $(\bu,q)\in H(\div;\Omega)\times
L^2(\Omega)$ such that
\begin{equation}
\begin{split}
(\beta\bu,\bv)+(q,\div\bv) & = (\bg,\bv)\quad\quad\forall\;\bv\in H(\div;\Omega),\\
(\div\bu,p) - (\alpha^{-1}q,p) & = 0\quad\quad\quad\quad\forall\;p\in L^2(\Omega).
\end{split}
\label{eq:mixed_weak}
\end{equation}
Notice that \eqref{eq:mixed_weak} coincides with the mixed formulation of the
scalar second order elliptic problem
\begin{equation}
-\div(\beta^{-1}\nabla q) + \alpha^{-1}q = g
\label{eq:pde_h1}
\end{equation}
by letting $\bu=\beta^{-1}\nabla q$. Hence, \eqref{eq:mixed_weak} can be
discretized by the hybridized mixed finite element methods designed for
\eqref{eq:pde_h1}, c.f. \cite{cg04}, resulting in a discrete problem of the form
\begin{equation*} \label{hyb3-ls}
\begin{bmatrix}
\hM & \hB^T & C^T \\
\hB & -\hW & 0 \\
C & 0 & 0
\end{bmatrix}
\begin{bmatrix}
\hx \\
q \\
\blambda
\end{bmatrix}
=
\begin{bmatrix}
\hf \\
0 \\
0
\end{bmatrix}.
\end{equation*}
The corresponding reduced problem for the Lagrange multipliers reads
\begin{equation*} \label{lag3-ls}
{\widetilde H}\blambda=
\begin{bmatrix}
C  & 0
\end{bmatrix}
\begin{bmatrix}
\hM & \hB^T \\
\hB & -\hW\\
\end{bmatrix}^{-1}
\begin{bmatrix}
C^T  \\ 0
\end{bmatrix}
\blambda =
\begin{bmatrix}
C  & 0
\end{bmatrix}
\begin{bmatrix}
\hM & \hB^T \\
\hB & -\hW\\
\end{bmatrix}^{-1}
\begin{bmatrix}
\hf \\ 0
\end{bmatrix},
\end{equation*}
which is identical to \eqref{Lagrange multiplier system} because of
\eqref{eq:discrete_hdiv} and the fact that the $(1,1)$ entry in the above
inverse is the inverse of the Schur complement $\hM + \hB^T \hW^{-1} \hB$.  In
other words, $H = {\widetilde H}$, so the algebraic hybridization approach is
identical when applied to the mixed or second order version of the problem.

\section{Scalable algebraic $H(\div)$ preconditioning} \label{section:preconditioning}

In this section we propose several scalable algebraic multigrid approaches for
finite element discretizations of the $H(\div)$ form $a_{\div}(\cdot,\cdot)$
based on the methods discussed in the previous sections.  Due to their
scalability on challenging practical problems, we base our preconditioners on
the parallel algebraic multigrid approaches available in the {\em hypre}
library, \cite{scaling_2012}.  Specifically we take advantage of the BoomerAMG
\cite{hy02}, AMS \cite{kolev-vassilevski-ams} and ADS \cite{kv12} methods in
{\em hypre}, targeting $H^1$, $H(\curl)$ and $H(\div)$ discretizations
respectively.

We first consider the static condensation Schur complement $S$, \eqref{static
condensation matrix}, and discuss approaches for its preconditioning.  Since
$S$ acts in general on traces of functions from a given space ($H(\div)$,
$H(\curl)$, or $H^1$), a suitable preconditioner for $S$ will be some
approximation to the Schur complement of an (optimal) preconditioner for the
original matrix $A$.  Of course, for efficiency, we would like to build a
preconditioner directly for the (sparse) matrix $S$ and not for the original
matrix $A$.  Such preconditioners are discussed in \cite{brunner-kolev-amg-explicit}
where it is shown, both theoretically and practically, that BoomerAMG and AMS
work well on Schur complements of $H^1$ and $H(\curl)$ discretizations
respectively. Similar analysis for the $H(\div)$ case was recently performed in
\cite{dpg-2016}.  In all three cases, the summary of the stable decomposition
theory is that if BoomerAMG/AMS/ADS works for a particular problem, it will also
work, purely algebraically, for a Schur complement of the problem.  Thus, the
$H(\div)$ Schur complement $S$ considered in this paper can be efficiently
preconditioned simply with a direct application of ADS.

We next consider preconditioning approaches for the hybridized matrix $H$.
Based on the connection with mixed methods from the previous section we argue
that one suitable preconditioner for the hybridized Schur complement $H$ in the
case of the $H(\div)$ problem is an $H^1$-based AMG built from the matrix $H$.
The rest of this section is devoted by providing arguments for this choice, by
exploring the ``near-nullspace'' of $H$, i.e. the set of its low-frequency
eigenmodes, that is critical in AMG theory.

Let us consider a norm $\vertiii{\cdot}$ for the space of Lagrange multiplier,
$$
{\vertiii{\blambda}}^2 =
  \sum_{\tau\in\mathcal{T}_h}\left(h_\tau\|\blambda\|^2_{\partial\tau}+
  \frac1{h_\tau}\|\blambda-m_\tau(\blambda)\|^2_{\partial\tau}\right),
$$
where $\mathcal{T}_h$ is a finite element mesh and $m_\tau(\blambda) =
\frac1{|\partial\tau|}\int_{\partial\tau}\blambda\;ds$.  It is proven in
\cite{cg05} that ${\widetilde H}$ is spectrally equivalent to $\vertiii{\cdot}$,
so we can deduce that $H$ has a small near-nullspace containing only the constant
functions, as opposed to the large near-nullspace of the original system $A$.
Therefore, the hybridized Schur complement $H$ is ``$H^1$-like'', which
motivates the use of $H^1$-based AMG for preconditioning.

In fact, by exploiting $\widetilde H$, one can see that $H$ has the explicit form:
$$
H = \widetilde H = C(\hM^{-1}-\hM^{-1}\hB^T(\hW+\hB\hM^{-1}\hB^T)^{-1}\hB\hM^{-1})C^T.
$$
In the above expression, $\hW$ is an $L^2$ mass matrix while $\hB\hM^{-1}\hB^T$
resembles a weighted discrete Laplacian. Hence, $\hB\hM^{-1}\hB^T \gg \hW$ and
so the near-nullspace of $H$ is characterized by the nullspace of the matrix
$C\hY C^T$, where
$$
\hY = \hM^{-1}-\hM^{-1}\hB^T(\hB\hM^{-1}\hB^T)^{-1}\hB\hM^{-1}.
$$
Note that for any $\hx\in \Ran(\hB^T)$, there exists $\hy$ such that $\hx=\hB^T\hy$, so
$$
\hY\hx = \hM^{-1}\hx-\hM^{-1}\hB^T(\hB\hM^{-1}\hB^T)^{-1}\hB\hM^{-1}\hB^T\hy =
  \hM^{-1}\hx-\hM^{-1}\hB^T\hy = 0.
$$
Thus, we have
\begin{equation}
\Ran(\hB^T)\subseteq \Ker(\hY).
\label{eq:inclusion1}
\end{equation}
On the other hand, let $\hx\in \Ker(\hY)$. Then
$$
\hM^{-1}\hx-\hM^{-1}\hB^T(\hB\hM^{-1}\hB^T)^{-1}\hB\hM^{-1}\hx = 0.
$$
Consequently, $\hx=\hB^T(\hB\hM^{-1}\hB^T)^{-1}\hB\hM^{-1}\hx.$ So for any
$\hy\in \Ker(\hB)$, we have
\begin{equation*}
\begin{split}
\left<\hx,\hy\right> & = \left<\hB^T(\hB\hM^{-1}\hB^T)^{-1}\hB\hM^{-1}\hx,\hy\right>\\
& = \left<(\hB\hM^{-1}\hB^T)^{-1}\hB\hM^{-1}\hx,\hB\hy\right> = 0
\end{split}
\end{equation*}
Hence,
$$
\Ker(\hY)\subseteq \Ker(\hB)^\perp = \Ran(\hB^T).
$$
The above inclusion relation and \eqref{eq:inclusion1} imply that
$$
\Ker(\hY) = \Ran(\hB^T).
$$
Therefore, in order to know what $\Ker(C\hY C^T)$ is, we should find out
$$
\Ran(C^T) \cap \Ran(\hB^T).
$$
Note that the constraint matrix $C$ enforces the continuity of the master and
slave dofs on the interface (cf.\ Example \ref{example: P and C}), so
$\Ran(C^T)$ is orthogonal to the subspace spanned by all the interior dofs. In
other words, $\hx$ is in $\Ran(C^T)$ only if every entry corresponding to an
interior dof is zero. So we are looking for the vectors in $\Ran(\hB^T)$ such
that entries corresponding to interior dofs are zero. Recall that $\hB$ comes
from the bilinear form $(\div\bu,p)$. In the case of Raviart-Thomas elements,
constants are the only $p$ such that $(\div\bu,p)=0$ for all the $\bu$
associated with interior dofs (bubbles). Hence, we deduce that
$$
\Ker(C\hY C^T) = \left\{\;\blambda\; |\; C^T\blambda\in\text{span}\{\hB^T\bone\}\;\right\},
$$
where $\bone$ is the coefficient vector of the constant $1$ function in the
corresponding discrete space. This clearly shows that
$\dim\big(\Ker(C\hY C^T)\big)\le 1$.

\section{Implementation details}\label{section: implementation}

In this section we provide some details for the practical application of the
presented approaches as implemented in the finite element library MFEM
\cite{mfem}.

The implementation is contained in class {\tt Hybridization} which can be used
either directly, or more conveniently through the {\tt BilinearForm} class.  The
prolongation operator $P$ is not directly constructed, instead, its action can
be performed using the methods
in class {\tt FiniteElementSpace}.  The constraints matrix $C$ is constructed as
illustrated in Example \ref{example: C from bilinear form}, based on a specified
multiplier finite element space and a given bilinear form.

\begin{figure}
\centering
\includegraphics[scale=0.75]{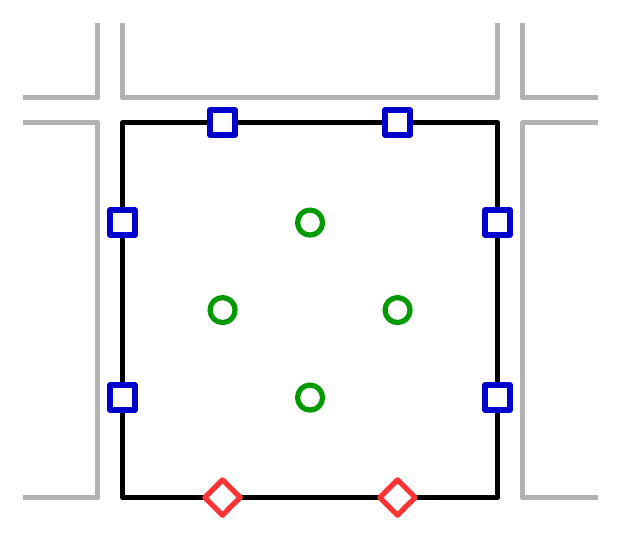}
\raisebox{0.4\height}{
\includegraphics[scale=0.75]{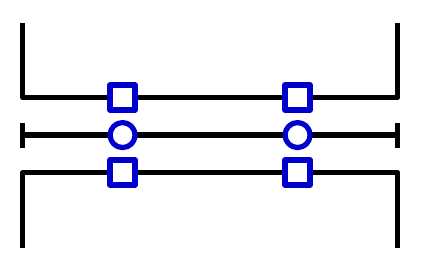}}
\raisebox{0.4\height}{
\includegraphics[scale=0.75]{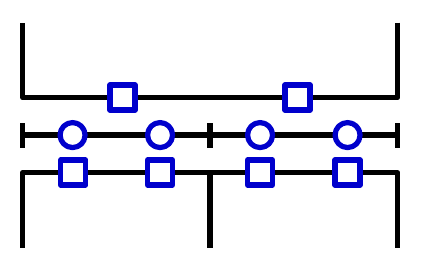}}
\caption{Left: local $RT_1$ dofs: $i$-dofs (green circles), $b$-dofs (blue
squares), and $e$-dofs (red diamonds). Center\slash Right: dofs associated with
a conforming (center) and a non-conforming (right) mesh face: $RT_1$ dofs to
match (blue squares) and multiplier dofs (blue circles).}
\label{fig:rt1-example}
\end{figure}

During a construction stage, the element matrices are reordered and stored in
the $2 \times 2$ block form \eqref{A-2x2}, see Figure \ref{fig:rt1-example} for
illustration. Recall that the $i$- and $b$-block unknowns correspond to the
element interior and boundary (or interface) degrees of freedom.  Any rows and
columns associated with essential degrees of freedom (e-dofs) are
eliminated. Algebraically, the interior dofs can be defined as those
corresponding to zero columns in the $C$ matrix. In other words, $C$ has the $1
\times 2$ block form
\[
C = \begin{bmatrix}  0 & C_{b} \end{bmatrix}.
\]
As a consequence, the matrix $H$ of the hybridized system
\eqref{Lagrange multiplier system} can be computed as
\[
H = C \wh{A}^{-1} C^T = C_b \wh{S}_{b}^{-1} C_{b}^T \,,
\]
where $\wh{S}_{b} = \wh{A}_{bb} - \wh{A}_{bi} \wh{A}_{ii}^{-1} \wh{A}_{ib}$ is the
Schur complement of $\wh{A}$ which can be formed independently for each
element in the mesh. In MFEM, this is performed using element-wise block LU
factorizations of the form
\[
\begin{bmatrix}
\wh{A}_{ii} & \wh{A}_{ib} \\
\wh{A}_{bi} & \wh{A}_{bb}
\end{bmatrix}
=
\begin{bmatrix}
\wh{L}_{i} & 0 \\
\wh{L}_{bi} & \wh{I}_{b}
\end{bmatrix}
\begin{bmatrix}
\wh{U}_{i} & \wh{U}_{ib} \\
 0 & \wh{S}_{b}
\end{bmatrix}
\]
to compute the Schur complements $\wh{S}_{b}$ which are then themselves
LU-factorized, allowing simple multiplication of any vector or matrix with
$\wh{S}_{b}^{-1}$.

Note that for Robin-type boundary conditions, one needs to assemble
contributions to the global system matrix A coming from integration over
boundary faces. In MFEM, such face matrices are assembled into the element
matrix of the element adjacent to the corresponding face. The rest of the
hybridization procedure remains the same.

\section{Numerical tests}\label{section: numerical tests}

To validate the parallel scalability of the proposed solvers, in this section we
present several weak scaling tests, where the number of degrees of freedom per
processor is held to be about the same while the number of processors is
increased. All the experiments except the ones in
Section~\ref{sec:radiation_diffusion} are performed on the cluster Sierra at
Lawrence Livermore National Laboratory (LLNL). Sierra has a total of 1944 nodes
(Intel Xeon EP X5660 @ 2.80 GHz), each of which has 12 cores and 24 GB of
memory. The numerical results are generated using the libraries MFEM \cite{mfem}
and {\em hypre} \cite{hypre} developed at LLNL.

\subsection{Soft hard materials}\label{sec:softhard}

Consider solving \eqref{eq:weak_hdiv} on the unit cube, i.e. $\Omega = [0,1]^3$,
with the source function $\bg$ being $[1, 1, 1]^T$. Let $\Omega_i =
[\frac14,\frac12]^3\cup[\frac12,\frac34]^3$. We consider $\alpha\equiv1$ in
$\Omega$, whereas
$$
\beta=\left\{
\begin{array}{ll}1 & \mbox{ in } \Omega\backslash\Omega_i\\10^p &
  \mbox{ in } \Omega_i\end{array}
\right., \;p \in\{ -8, -4, 0, 4, 8\}.
$$
\begin{figure}[h!]
\centering
\includegraphics[scale=.14,clip,trim=0 0 0 1cm]{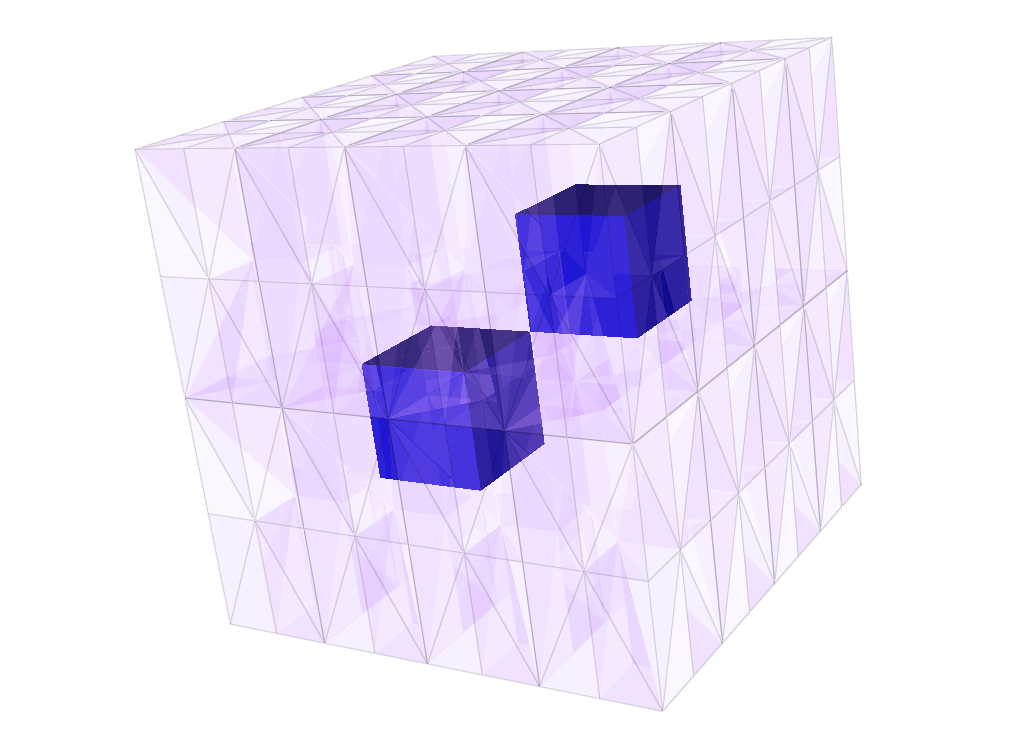}
\end{figure}
This example is intended to test the robustness of the solvers against
coefficient jumps. We discretize \eqref{eq:weak_hdiv} by Raviart-Thomas elements
($RT_k$) on hexahedral meshes.  In order to better contrast the effect of higher
order discretizations, the initial meshes in the weak scaling tests are chosen
so that the problem size is the same for discretizations of different orders,
see Table~\ref{tab:order}.
\begin{table}[h!]
\centering
\caption{Initial mesh sizes in the weak scaling tests in Section~\ref{sec:softhard}}
\label{tab:order}
\begin{tabular}{|c||c|c|c|} \hline
Finite element & $RT_0 $ & $RT_1 $ & $RT_3 $ \\
 \hline
Initial mesh &   $64\times64\times32$ & $32\times32\times16$ & $16\times16\times8$ \\
\hline
\end{tabular}
\end{table}
The PCG stopping criterion is a reduction in the $l_2$ residual by a factor of
$10^{-12}$.  The time to solution and number of PCG iterations (in parentheses)
for the proposed approaches, as well as ADS, are reported in
Tables~\ref{tab:softhard_rt0}--\ref{tab:softhard_rt3}. Here, the time to
solution for the hybridization approach includes the formation time of the
hybridization Schur complement $H$, the setup time of AMG (BoomerAMG
\cite{hy02}), PCG solving time, and the time for recovering the original
solution by back substitution.  Similarly, the time to solution for the static
condensation approach includes the formation time of $S$
(cf. Section~\ref{section: static condensation}), ADS setup time, PCG solving
time, and back substitution. Lastly, the time to solution for ``ADS'' is simply
ADS setup time and PCG solving time.
\begin{table}[h!]
\centering
\caption{Time to solution of different solvers: $RT_0$ on hexahedral meshes}
\label{tab:softhard_rt0}
\scalebox{0.9}{
\begin{tabular}{ccccccc} \hline
\#proc & $\;$Size$\;$ & $\;p=-8\;$ & $\;p=-4\;$ &  $\;p=0\;$ &  $\;p=4\;$ & $\;p=8\;$\\
 \hline
\multicolumn{7}{c}{Hybridization + AMG}\\
\hline
3 &	401,408 &	3.39	(24)	&	3.5	(26)	&	3.55	(27)	&	3.56	(27)	&	3.44	(26)	\\
24 &	3,178,496 &	7.29	(29)	&	7.43	(31)	&	7.34	(30)	&	7.29	(30)	&	7.35	(30)	\\
192 &	25,296,896 &	8.28	(33)	&	8.54	(37)	&	8.45	(32)	&	8.3	(32)	&	8.3	(33)	\\
1,536 &	201,850,880 &	9.79	(37)	&	10.33	(39)	&	9.98	(36)	&	10.28	(38)	&	10.17	(36)	\\
 \hline
\multicolumn{7}{c}{ADS}\\
\hline
3 &	401,408 &	17.35	(16)	&	17.51	(16)	&	15.94	(13)	&	20.54	(19)	&	20.11	(21)	\\
24 &	3,178,496 &	36.94	(18)	&	36.91	(18)	&	35.69	(16)	&	38.76	(21)	&	40	(23)	\\
192 &	25,296,896 &	42.54	(20)	&	42.39	(20)	&	40.99	(17)	&	43.75	(22)	&	45.7	(25)	\\
1,536 &	201,850,880 &	49.91	(21)	&	49.64	(21)	&	47.73	(19)	&	52.93	(25)	&	52.62	(27)	\\
\hline
\end{tabular}}
\end{table}
\begin{table}[h!]
\centering
\caption{Time to solution of different solvers: $RT_1$ on hexahedral meshes}
\label{tab:softhard_rt1}
\scalebox{0.9}{
\begin{tabular}{ccccccc} \hline
\#proc & $\;$Size$\;$ & $\;p=-8\;$ & $\;p=-4\;$ &  $\;p=0\;$ &  $\;p=4\;$ & $\;p=8\;$\\
 \hline
\multicolumn{7}{c}{Hybridization + AMG}\\
\hline
3 &	401,408 &	3.35	(27)	&	3.39	(28)	&	3.46	(28)	&	3.43	(29)	&	3.46	(28)	\\
24 &	3,178,496 &	5.97	(30)	&	6.14	(33)	&	6.06	(32)	&	5.95	(31)	&	6.11	(32)	\\
192 &	25,296,896 &	6.85	(36)	&	7.03	(37)	&	7.25	(36)	&	7	(35)	&	7.43	(35)	\\
1,536 &	201,850,880 &	7.84	(39)	&	9.17	(41)	&	8.17	(39)	&	8.83	(41)	&	9.85	(40)	\\
\hline
\multicolumn{7}{c}{Static condensation + ADS}\\
\hline
3 &	401,408 &	19.69	(15)	&	19.65	(15)	&	18.23	(12)	&	20.92	(17)	&	21.91	(18)	\\
24 &	3,178,496 &	40.5	(19)	&	40.69	(19)	&	36.77	(15)	&	42.42	(21)	&	43.51	(22)	\\
192 &	25,296,896 &	46.48	(21)	&	46.81	(21)	&	42.63	(17)	&	48.68	(23)	&	51.62	(26)	\\
1,536 &	201,850,880 &	54.78	(23)	&	54.43	(22)	&	50.38	(19)	&	59.4	(26)	&	61.56	(29)	\\
\hline
\multicolumn{7}{c}{ADS}\\
\hline
3 &	401,408 &	34.73	(18)	&	34.53	(17)	&	30.66	(14)	&	38.93	(22)	&	41.44	(24)	\\
24 &	3,178,496 &	60.9	(20)	&	61.13	(20)	&	55.19	(16)	&	68.06	(25)	&	71.44	(27)	\\
192 &	25,296,896 &	68.55	(22)	&	67.39	(21)	&	62.5	(18)	&	77.74	(28)	&	81.1	(30)	\\
1,536 &	201,850,880 &	78.8	(24)	&	80.99	(24)	&	72.89	(20)	&	90.74	(31)	&	102.93	(37)	\\
 \hline
\end{tabular}}
\end{table}
\begin{table}[h!]
\centering
\caption{Time to solution of different solvers: $RT_3$ on hexahedral meshes}
\label{tab:softhard_rt3}
\scalebox{0.87}{
\begin{tabular}{ccccccc} \hline
\#proc & $\;$Size$\;$ & $\;p=-8\;$ & $\;p=-4\;$ &  $\;p=0\;$ &  $\;p=4\;$ & $\;p=8\;$\\
 \hline
\multicolumn{7}{c}{Hybridization + AMG}\\
\hline
3 &	401,408 &	9.27	(33)	&	9.46	(35)	&	9.42	(37)	&	9.62	(38)	&	9.51	(38)	\\
24 &	3,178,496 &	13.1	(41)	&	12.91	(43)	&	14	(44)	&	13.53	(44)	&	13.25	(44)	\\
192 &	25,296,896 &	14.81	(48)	&	15.58	(50)	&	15.26	(50)	&	14.71	(50)	&	16.4	(50)	\\
1,536 &	201,850,880 &	16.8	(53)	&	16.33	(56)	&	16.63	(56)	&	17.67	(56)	&	19.06	(56)	\\
\hline
\multicolumn{7}{c}{Static condensation + ADS}\\
\hline
3 &	401,408 &	31.38	(12)	&	31.38	(12)	&	29.28	(10)	&	32.65	(13)	&	32.38	(13)	\\
24 &	3,178,496 &	56.09	(15)	&	55.87	(15)	&	54.64	(14)	&	59.06	(17)	&	59.03	(17)	\\
192 &	25,296,896 &	67.01	(18)	&	67.28	(18)	&	65.22	(17)	&	70.43	(20)	&	70.53	(20)	\\
1,536 &	201,850,880 &	83.81	(20)	&	82.14	(20)	&	81.69	(19)	&	85.65	(22)	&	86.65	(23)	\\
\hline
\multicolumn{7}{c}{ADS}\\
\hline
3 &	401,408 &	184.67	(18)	&	186.1	(18)	&	169.27	(14)	&	188.58	(19)	&	189.52	(20)	\\
24 &	3,178,496 &	253.3	(19)	&	259.11	(20)	&	240.87	(16)	&	269.12	(22)	&	284.04	(25)	\\
192 &	25,296,896 &	288.32	(22)	&	289.14	(22)	&	261.95	(17)	&	305.86	(25)	&	325.93	(29)	\\
1,536 &	201,850,880 &	329.76	(24)	&	330.46	(23)	&	301.61	(19)	&	353.63	(28)	&	366.94	(31)	\\
 \hline
\end{tabular}}
\end{table}
Our results show that all the solvers have good weak scaling and robustness with
respect to coefficient jumps. In all the cases, while both the hybridization and
static condensation solvers give shorter solving time than ADS, the
hybridization approach is the fastest. Moreover, as the finite element order
goes higher, discrepancy in the solution time between the solvers becomes more
significant even though the problem size is the same. In particular, for the
$RT_3$ discretization, hybridization is in general 4 times faster than static
condensation, and 20 times faster than ADS.
\begin{figure}[h!]
\centering
\includegraphics[scale=.3,clip,trim=3.8cm 0 4.8cm 0]{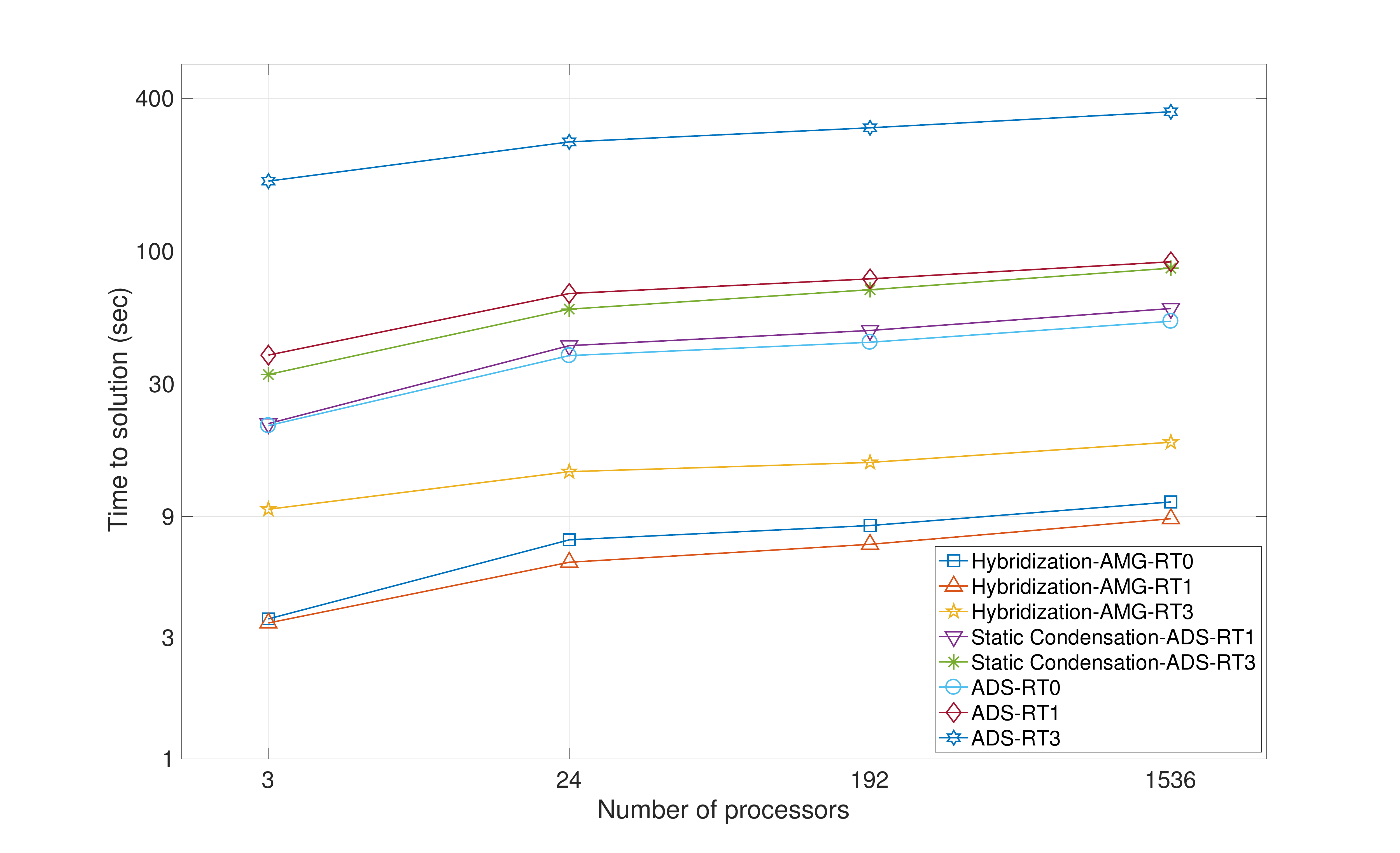}
\caption{Weak scaling plots in $\log$-$\log$ scale for the soft hard materials
  example when $p=4$}
\label{fig:weakscaling_all}
\end{figure}

\subsection{The crooked pipe problem}

\begin{figure}[h!]
\centering
\subfigure{\includegraphics[scale=.11]{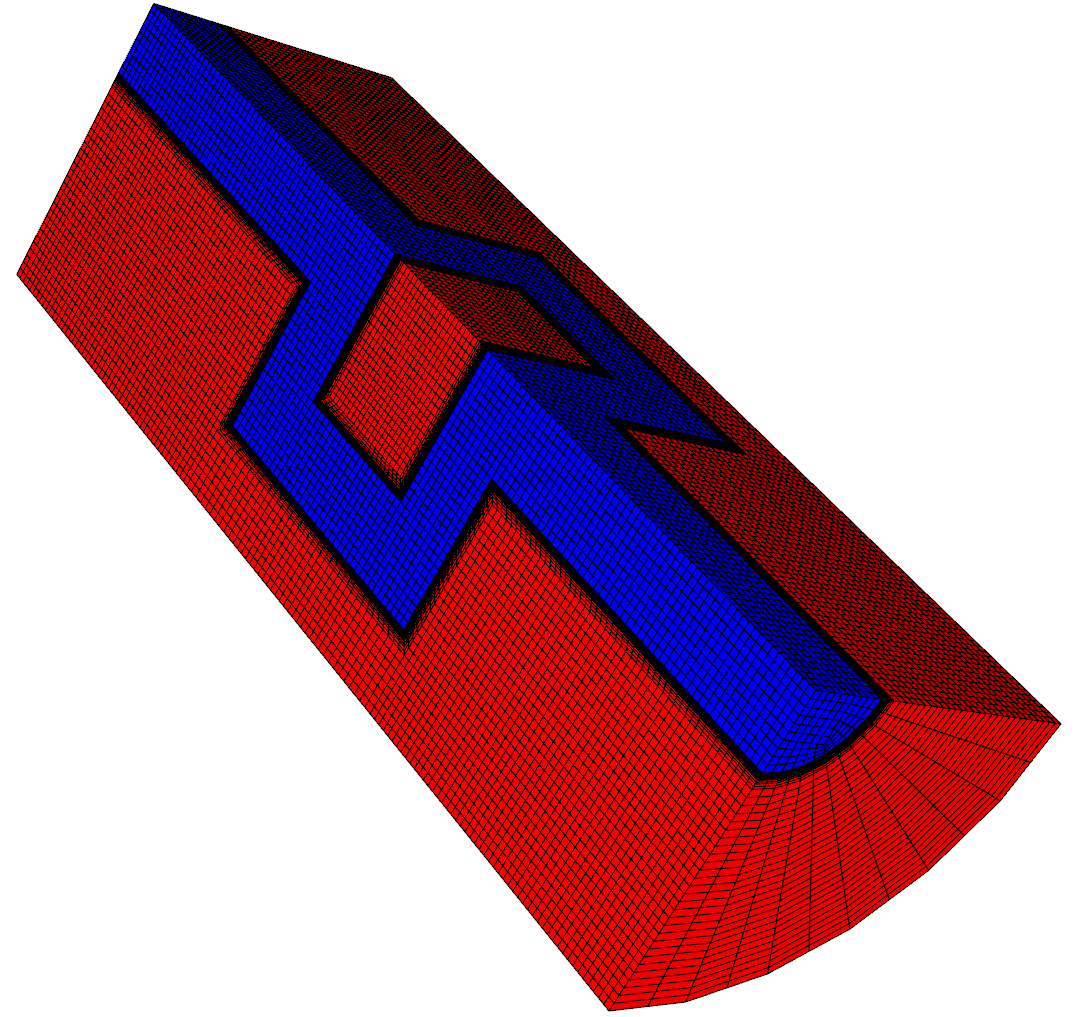}}\hspace{15mm}
\subfigure{\includegraphics[scale=.13,clip,trim=0 3cm 0 0]{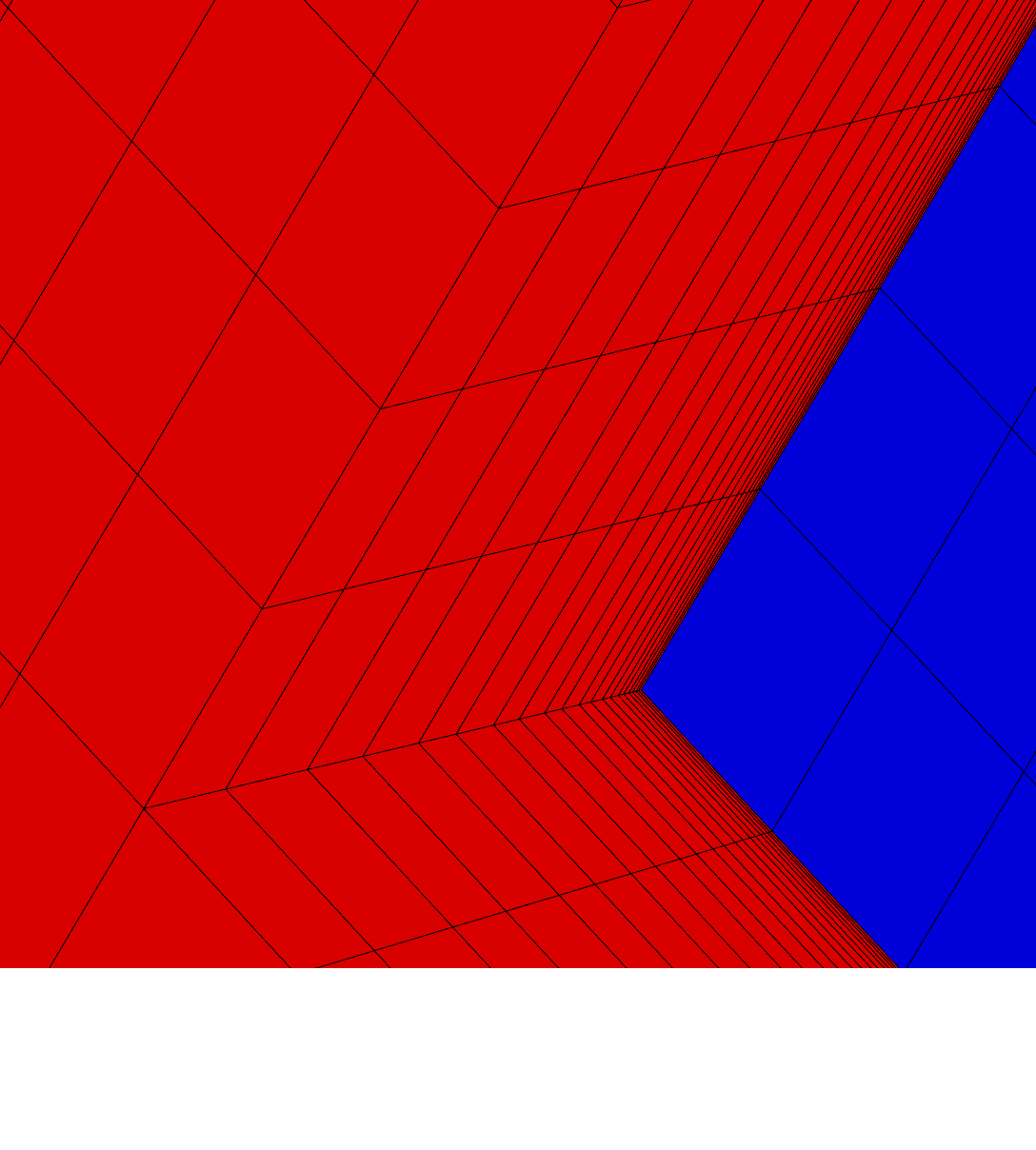}}
\caption{Geometry and material interface of the crooked pipe example}
\label{fig:crookedpipe}
\end{figure}

This example deals with large coefficient jumps and highly stretched elements,
see Figure~\ref{fig:crookedpipe}. The coefficients $\alpha$ and $\beta$ are both
piecewise-constant functions such that
\[
(\alpha, \beta)=\left\{
\begin{array}{ll}(1.641,0.2) & \mbox{ in red region,} \\
(0.00188,2000) & \mbox{ in blue region.}  \end{array}
\right.
\]
The setup of the experiments is the same as in Section~\ref{sec:softhard}, and
the weak scaling results are reported in Table~\ref{tab:crookedpipe}. Note that
we were running out of memory when setting up ADS for $RT_3$ (indicated by
``$\ast$"), which shows that the hybridization and static condensation
approaches are more memory-friendly for high order discretizations.  The results
in Table~\ref{tab:crookedpipe} again suggest that hybridization with AMG is the
most efficient solver among the three for this particular problem; it is also
observed that the iteration count in the hybridization approach is more steady
than the other two.

\begin{table}[h!]
\centering
\caption{Solution time for the crooked pipe example}\label{tab:crookedpipe}
\scalebox{0.85}{
\begin{tabular}{ccccc} \hline
 \#proc$\;$ & $\;$Size$\;$ & $\;$\quad\quad$RT_0$\quad\quad$\;$ &
   $\;$\quad\quad$RT_1$\quad\quad$\;$ &  $\;$\quad\quad$RT_3$\quad\quad$\;$ \\
 \hline
\multicolumn{5}{c}{Hybridization + AMG}\\
\hline
12 &	2,805,520 &	8.56	(34)	&	8.44	(33)	&	22.41	(42)	\\
96 &	22,258,240 &	13.4	(37)	&	12.45	(36)	&	29.97	(44)	\\
768 &	177,322,240 &	17.82	(45)	&	15.28	(45)	&	37.13	(49)	\\
6,144 &	1,415,603,200 &	23.73	(54)	&	21.36	(55)	&	43.01	(60)	\\
\hline
\multicolumn{5}{c}{Static condensation + ADS}\\
\hline
12 &	2,805,520 &	73.76	(52)	&	98.73	(43)	&	190.53	(40)	\\
96 &	22,258,240 &	126.87	(80)	&	177.71	(74)	&	305.51	(67)	\\
768 &	177,322,240 &	191.67	(117)	&	291.03	(122)	&	481.28	(107)	\\
6,144 &	1,415,603,200 &	274.09	(168)	&	463.95	(190)	&	715.79	(158)	\\
\hline
\multicolumn{5}{c}{ADS}\\
\hline
12 &	2,805,520 &	73.76	(52)	&	174.55	(59)	&	$\ast$	($\ast$)	\\
96 &	22,258,240 &	126.87	(80)	&	291.4	(93)	&	$\ast$	($\ast$)	\\
768 &	177,322,240 &	191.67	(117)	&	442	(139)	&	$\ast$	($\ast$)	\\
6,144 &	1,415,603,200 &	274.09	(168)	&	706.5	(212)	&	$\ast$	($\ast$)	\\
 \hline
\end{tabular}}
\end{table}

\subsection{Application to radiation diffusion}\label{sec:radiation_diffusion}

\begin{figure}[h!]
\centering
\vspace{5mm}
\includegraphics[scale=.23,clip, trim = 3cm 3cm 1cm 3cm]{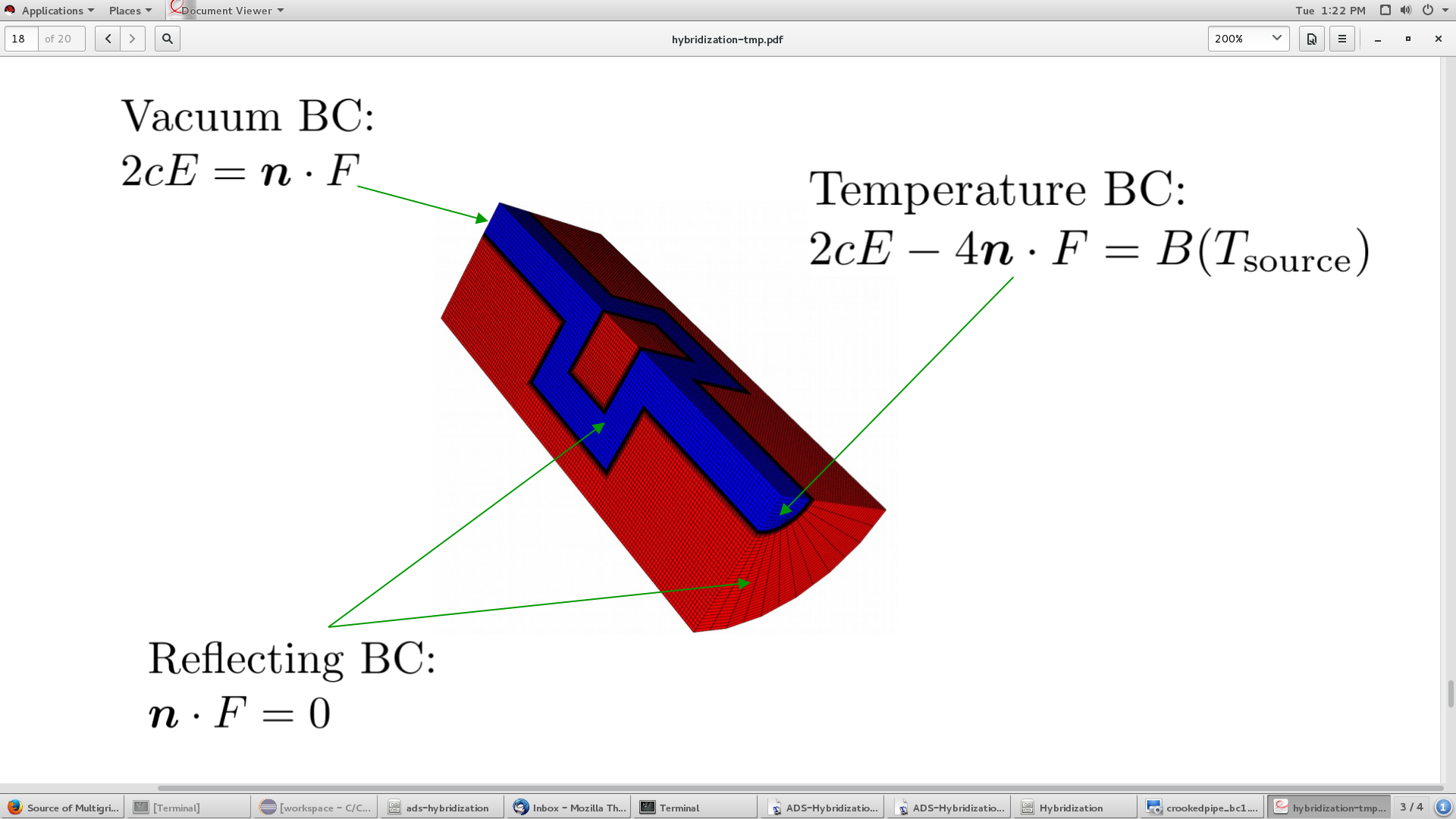}
\caption{Boundary conditions for the radiation diffusion problem}
\label{fig:crookedpipe_bc}
\end{figure}

In this section we compare the proposed approaches when they are used for
solving the radiation diffusion equations in flux form, see \cite{Brunner02}.
This model is often used in astrophysical and inertial confinement fusion (ICF)
simulations. We consider a problem with two materials ($k=1,2$ below), indicated
by blue and red colors in Figure~\ref{fig:crookedpipe_bc}, and there are no
mixed zones.  The governing equations are
\begin{equation}
\begin{split}
  \rho_k \frac{de_k}{dt} &= c \sigma_{k} (E-B(T(e_k))), \\
  \frac{dE}{dt}+\nabla \cdot F & = -c \sum_k \sigma_{k}(E-B(T(e_k))), \\
  \frac13 \nabla E &= -\sum_k \frac{\sigma_{k}} c F,
\end{split}
\end{equation}
subject to the boundary conditions specified in Figure~\ref{fig:crookedpipe_bc}.
Here $e_k, \rho_k, \sigma_k$ are the internal energy, density and opacity,
respectively, of material $k$, $E$ is radiation energy, $F$ is radiation flux,
$B$ is the Planck's function for black-body radiation, and $c$ is the speed of
light.  The temperature boundary condition heats up one of the ends of the pipe,
causing an increase of the radiation field $E$ and its diffusion throughout the
domain.

The discrete approximations of $E$ and $e_k$ are in the finite elements space
$Q_r \subset L^2$, consisting of piecewise polynomials of degree $r$, while the
flux $F$ is approximated in the Raviart-Thomas space $RT_r$.  For stability
reasons the time stepping of $e_k$ and $E$ is implicit, which results in a
nonlinear problem in each time step, as $B(T(e_k))$ is nonlinear.  Newton's
method is used to solve the resulting system. As $E$ and $e_k$ are
discontinuous, they can be eliminated locally during each Newton iteration,
which is inexpensive and inherently parallel. Upon eliminating $E$ and $e_k$, an
$H(\div)$ linear system for $F$ is obtained, which we use to test the
performance of the two proposed approaches, as well as AMS (2D) and ADS (3D).


Because each method ultimately solves a different linear system, this example is
used to obtain better perception of the execution time needed for each method to
obtain similar solution accuracy with respect to the original system.  Starting
with identical initial conditions, each computation performs exactly one linear
solve, and then the $F$ residual is calculated from the definition of the
nonlinear system. Execution times are presented for computations reaching
similar residuals. For this particular example, this is achieved by setting the
relative tolerances to 1e-6, 1e-5 and 1e-12 for AMS/ADS, static
condensation$+$AMS/ADS, and hybridization$+$AMG, respectively.

As observed in Tables~\ref{tab:rd2D} and \ref{tab:rd3D}, even though the
hybridized system must be solved to a much lower relative tolerance, the
hybridization approach still has the shortest solving time, followed by the
static condensation and then ADS.  Except for $RT_0$ in the 2D case, the
hybridization approach is at least 10 times faster than AMS/ADS, which is a
substantial improvement.
\begin{table}[h!]
\centering
\caption{2D problem with \#dof=324712.
}
\scalebox{0.85}{
\begin{tabular}{|c||cc||ccc||ccc|} \hline
& \multicolumn{2}{|c||}{$RT_0$} & \multicolumn{3}{|c||}{$RT_1$} &
  \multicolumn{3}{c|}{$RT_3$}\\
\hline
Method & HB   & AMS   & HB   & SC   & AMS   & HB   & SC   & AMS   \\
 \hline
Time   & 3.05     & 15.35 & 1.11 & 7.29 & 18.63 & 1.35 & 6.21 & 35.68 \\
\giter & 100      & 245   & 50   & 129  & 200   & 73   & 130  & 190   \\
\fires & 2.09     & 5.34  & 9.65 & 3.76 & 6.64  & 3.61 & 3.58 & 2.01  \\
\hline
\end{tabular}}
\label{tab:rd2D}
\end{table}
\begin{table}[h!]
\centering
\caption{3D problem with \#dof=124208.}
\scalebox{0.85}{
\begin{tabular}{|c||cc||ccc||ccc|} \hline
& \multicolumn{2}{|c||}{$RT_0$} & \multicolumn{3}{|c||}{$RT_1$} &
  \multicolumn{3}{c|}{$RT_3$}\\
\hline
Method & HB   & ADS  & HB   & SC   & ADS   & HB   & SC   & ADS   \\
 \hline
Time   & 0.39 & 3.86 & 0.58 & 4.39 & 10.06 & 3.02 & 12.44 & 67.68 \\
\giter & 24   & 42   & 27   & 23   & 43    & 28   & 18    & 46    \\
\fires & 8.23 & 3.87 & 1.48 & 7.00 & 2.54  & 1.01 & 0.882 & 1.19  \\
\hline
\end{tabular}}
\label{tab:rd3D}
\end{table}

\subsection{The SPE10 benchmark}

\begin{figure}[h!]
\centering
\includegraphics[scale=.3,clip, trim = 0 8.5cm 5cm 10cm]{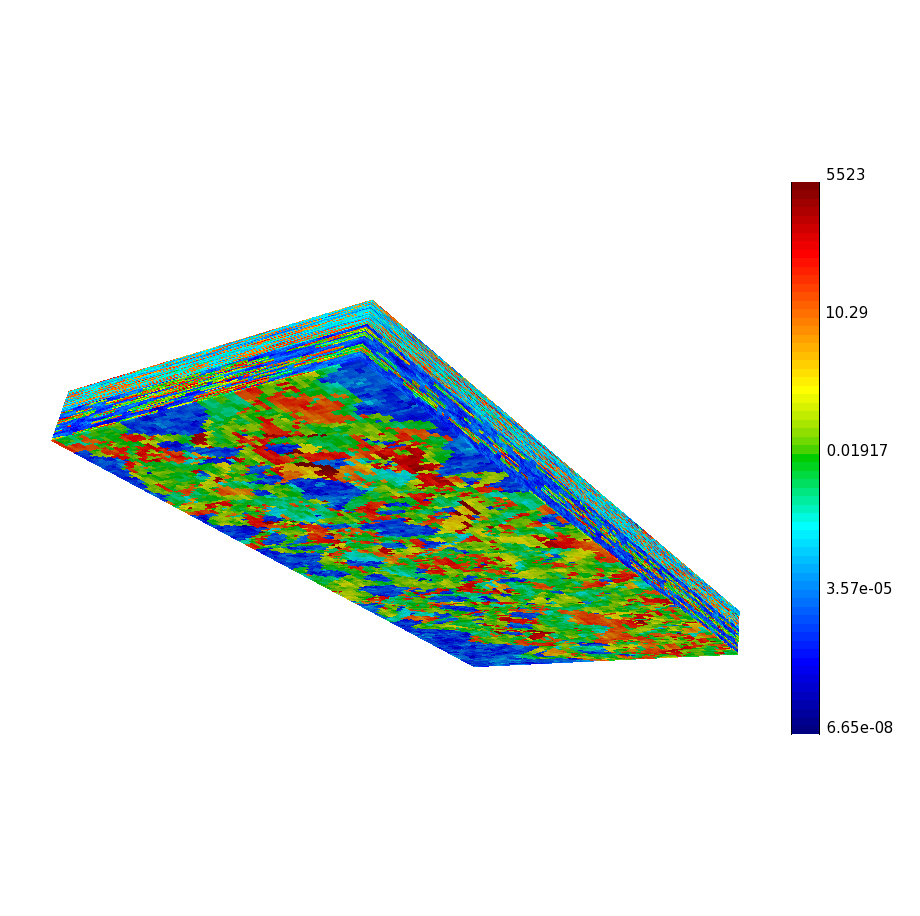}
\hspace{3mm}
\includegraphics[scale=.25,clip, trim = 20cm 3.8cm 0 3.6cm]{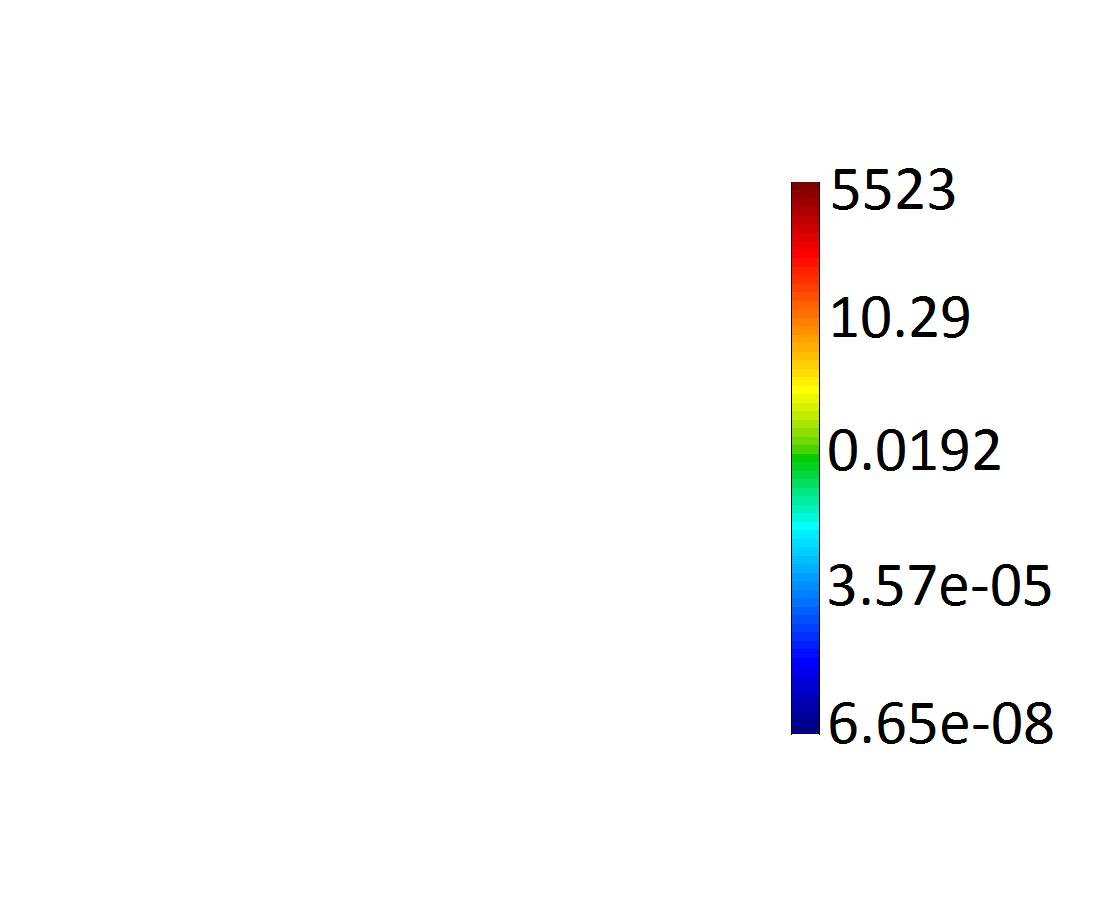}
\caption{The SPE10 coefficient}
\label{fig:spe10}
\end{figure}

Our last numerical example is the SPE10 benchmark coefficient \cite{spe10} from
the reservoir simulation community. This coefficient is well-known for its
multi-scale heterogeneity and anisotropy. The initial mesh is a
$60\times220\times85$ hexahedral grid, with the size of each cell being
$20\times10\times2$. $\beta$ is taken to be the SPE10 coefficient (see
Figure~\ref{fig:spe10}), and $\alpha \equiv1$. The time spent on different
components of the solution process for the hybridization approach and ADS are
shown in Tables~\ref{tab:spe10_hb} and \ref{tab:spe10_ads} respectively.  This
example also demonstrates the superiority in terms of performance of the
hybridization solver.
\begin{table}[h!]
\centering
\caption{The SPE10 example solved by hybridization and AMG}\label{tab:spe10_hb}
\scalebox{0.85}{
\begin{tabular}{ccccccc} \hline
 \#proc$\,$ & $\,$Size$\,$ & $\,$Hybridize$\,$  &  $\,$AMG setup$\,$  &
   $\,$PCG Solve$\,$  &  $\,$Back sub.$\,$ &  $\,$Total  time$\,$\\
\hline
3 &	3,403,000 &		21.3	&	2.21	&	18.6	&	0.51	&	42.62	(48)	\\
24 &	27,076,000 &	68.75	&	3.76	&	23.59	&	0.62	&	96.72	(45)	\\
192 &	216,016,000 &	111.82	&	6.93	&	27.85	&	0.63	&	147.22	(44)	\\
1,536 &	1,725,760,000 &	131.85	&	9.6	&	29.58	&	0.65	&	171.68	(40)	\\
\hline
\end{tabular}}
\end{table}
\begin{table}[h!]
\centering
\caption{The SPE10 example solved by ADS}\label{tab:spe10_ads}
\scalebox{0.85}{
\begin{tabular}{ccccc} \hline
 \#proc$\;$ & $\;$Size$\;$ & $\;$ADS setup$\;$ & $\;$PCG solve$\;$ &
   $\;$\quad Total time\quad$\;$\\
\hline
3 &	3,403,000 &	85.77	& 405.76	&	491.53	(123)\\
24 &	27,076,000 &	290.17	& 1416.49	&	1706.66	(300)	\\
192 &	216,016,000 &	492.11	& 3529.19	&	4021.3	(635)\\
1,536 &	1,725,760,000 &	585.7	& 5977.23	&	6562.93	(959)	\\
\hline
\end{tabular}}
\end{table}

\section{Conclusion} \label{section: conclusion}

In this paper, two popular techniques in finite elements, hybridization and
static condensation, were studied and compared. Both of these two techniques
lead to a reduced system of the same size. The formation of the reduced system
involves inversion of independent local matrices, which can be done in
parallel. We also discuss suitable preconditioning strategies using various
algebraic multigrid methods for the respective reduced systems. In particular,
we demonstrate that if BoomerAMG/AMS/ADS work well for the original problem,
then they also work well for the reduced systems obtained from static
condensation. On the other hand, AMG for $H^1$-equivalent problems is actually a
well-suited preconditioner for the reduced system resulting from hybridization
of $H(\div)$ problems. We presented several numerical examples comparing the
performance of the proposed approaches, as well as solving the original problem
directly by ADS/AMS. Our results show good weak scaling of all of the proposed
approaches. In general, for $H(\div)$ problems, the approach using hybridization
and AMG is more efficient than the approach using static condensation and
ADS/AMS, while both of these approaches are faster than a direct application of
ADS/AMS to the original problem. In fact, we observed substantial savings in the
solve time if the hybridization approach is used, especially for higher order
discretizations. Implementation of the two approaches are freely available in
MFEM \cite{mfem}.

\bibliographystyle{siam}    


\bibliography{references}       

\end{document}